\documentstyle[12pt]{article}
\def\baselinestretch{1.2}
\title{\bf Fuzzy $L$ languages
}
\author{  Zu-Guo Yu$^{1,2}$\\
 {\small $^{1}$Institute of Theoretical Physics, Academia Sinica},\\
  {\small P.O. Box 2735, Beijing 100080, China. Email: yuzg@itp.ac.cn}  \\
 {\small $^{2}$Department of Mathematics, Xiangtan University, Hunan 411105, China}.
  }
\date{}
\begin{document}
\newtheorem{Theorem}{\quad Theorem}
\newtheorem{Proposition}{\quad Proposition}[section]
\newtheorem{Definition}{\quad Definition}[section]
\newtheorem{Lemma}{\quad Lemma}[section]
\newtheorem{Corollary}{\quad Corollary}[section]
\newtheorem{Example}{\quad Example}[section]
\maketitle
\begin{abstract}
  A fuzzy aspect is introduced into $L$ systems and some of its properties
  are investigated. The relationship between fuzzy $L$ languages and the
  fuzzy languages generated by fuzzy grammar defined by Lee and Zadeh is
  elucidated. The concept of fuzzy entropy of a string is proposed as a measure
  of fuzziness of the latter. Some relationship between fuzzy $L$ systems and
  the ordinary $L$ systems is discussed.
\end{abstract}

\section{Introduction}
 
 {~~~~}Natural languages such as English have inexactitude and ambiguity both
 syntactically and semantically. A way of extending the concepts of formal
 languages to incorporate some aspects of natural language is the introduction
 of fuzziness. Lee and Zadeh$^{[3]}$, M. Mizumoto {\it et al} $^{[5]}$
 introduced fuzzy grammars as an extension of ordinary formal grammars (see
 e.g., Ref.[7]) by using the concept of fuzzy set$^{[11]}$. Their formal grammar
 is based on the Chomsky hierarchy$^{[1]}$. However, Chomsky hierarchy is not
 the only way to generate languages and to classify their complexity.
 There exists, for example, another grammatical hierarchy, called $L$-systems
 or developmental systems $^{[9]}$. The $L$-systems were first introduced by
 the biologist Lindenmayer$^{[4]}$ to model the growth and development of
 multicellular filamentous organisms. The $L$ systems are based on parallel
 rewriting rules. When a string is rewritten each symbol in the string must be
 renewed at the same time according to the rewriting rules. The standard
 references on $L$-systems are Herman and Rozenberg$^{[2]}$, Rozenberg and
 Salomaa$^{[6]}$.

  In this paper we introduce some families of fuzzy $L$-systems
and investigate their properties. We further discuss the relationship
between fuzzy $L$ languages and the fuzzy languages generated by fuzzy grammar
proposed in Ref.[3,5]. A measure of fuzziness for a string, called the
fuzzy entropy of a string with respect to a given fuzzy $L$ system, will be
defined. The relationship between fuzzy $L$ languages and the ordinary $L$
languages is also discussed.

  We assume the reader to be familiar with the rudiments of formal language
  theory (see, e.g., Salomaa$^{[8]}$) and of $L$ system theory (see, e.g.,
  Rozenberg and Salomaa$^{[6]}$). The reader may consult the book (Ref.[12])
  of Zimmermann for basic concepts of fuzzy set.

\section{Basic definitions and notations}

{~~~~}In this section we give some basic definitions on $L$ systems and
fix the notations to be used later on.

  An {\it alphabet} is a set of abstract symbols. The alphabets we consider are
always finite nonempty sets. The elements of an alphabet $\Sigma$ are called
{\it letters} or {\it symbols}. A {\it word} over an alphabet $\Sigma$ is a
finite string consisting of zero or more letters of $\Sigma$. The string
consisting of zero letters is called an empty word and is denoted as $\Lambda$.
 The set of all words
(resp. all nonempty words) over an alphabet $\Sigma$ is denoted by $\Sigma^*$
(resp. $\Sigma^+$). Subsets of $\Sigma^*$ are referred to as {\it languages} over
$\Sigma$. For words $w_1$ and $w_2$, the juxtaposition $w_1w_2$ is called the
{\it catenation} of $w_1$ and $w_2$. The {\it catenation} of two languages $L_1$
and $L_2$ is defined by
$$L_1L_2=\{w_1w_2|w_1\in L_1, w_2\in L_2\}.$$

  For each letter $a$ of an alphabet $\Sigma$, let $h(a)$ be a language, possibly
over a different alphabet. Define, further more, $$h(\Lambda)=\{\Lambda\},\quad
h(w_1w_2)=h(w_1)h(w_2).$$ For a language $L$ over $\Sigma$, we define
$$h(L)=\{u|u\in h(w)\ \mbox{for some}\ w\in L\}.$$
Such a mapping $h$ is called a {\it substitution}.

  A substitution $h$ such that each $h(a)$ consists of a single word is called
 a {\it homomorphism} or, briefly, a {\it morphism}.

  From Rozenberg and Salomaa [6], the finite substitution or homomorphism $h$
 will often be defined by listing the productions for each letter in $\Sigma$.
 Then we can view $h$ as a set of productions for each letter in $\Sigma$.
 In this paper we do not distinguish whether $h$ is a mapping or a set of
 productions for each letter in $\Sigma$.

\section{Fuzzy $L$ systems}
{~~~~}In this section we formulate a few elementary definitions and some
results of inclusion.

\bigskip   {\bf Definition 3.1}: A {\it fuzzy $0L$ system} $FG$ is a system
 $$FG=(\Sigma,\widetilde{h},\omega,J,f)$$
 where

  (i) $\Sigma$ is an alphabet.

  (ii) $\omega$, referred to as the {\it axiom}, is an element of $\Sigma^*$.

  (iii) $h$, viewed as the set of productions for each letter in $\Sigma$,
   is a finite substitution on $\Sigma$ (into the set of subsets of $\Sigma^*$);
   $\widetilde{h}$ is a set of productions defined by
  $$\widetilde{h}=\{(r)a\rightarrow uf(r)|,\quad a\in \Sigma\}$$
  where $r\in J$, $a\rightarrow u$ is an ordinary rewriting rule in $h$, $f(r)$ is
  the grade of the application of the production $r$, which will be defined
  in (v) below.

  (iv) $J$ is the set of production labels as shown in (iii), $J=\{r\}$.

  (v) $f$ is a membership function such as
  $$f:J\rightarrow [0,1].$$
  $f$ may be called a {\it fuzzy function}, and the value $f(r)$, $r\in J$, is
  the grade of the application of a production $r$.

   $(r)a\rightarrow uf(r)$ is called a
{\it  fuzzy production}. $\widetilde{h}$ is called {\it fuzzy substitution}
corresponding to substitution $h$ on
$\Sigma$.
   When $h$ is a homomorphism,
  $FG$ is called a {\it fuzzy $D0L$ system}.

\bigskip{\bf Example 3.1}: Let $FG=(\{a,b\},\widetilde{h},aba, J,f)$, the productions are
$$(1)\quad a\rightarrow aba\quad 0.5$$
$$(2)\quad b\rightarrow \Lambda\qquad 0.7,$$
then $FG$ is a fuzzy $D0L$ system.

{\bf Example 3.2}: Let $FG=(\{0\},\widetilde{h},0,\{(1),(2)\},f)$, where
$$\widetilde{h}=\{(1)\quad 0\rightarrow \Lambda\quad 0.3,\quad (2)\quad 0\rightarrow 00\quad 0.8\},$$
then $FG$ is a fuzzy $0L$ system.

\bigskip  {\bf Definition 3.2}: Let $FG=(\Sigma, \widetilde{h},\omega,J,f)$ be
 a fuzzy $0L$ system. A {\it Fuzzy derivation} $FD$ is a triple
 $({\cal O},\nu,p)$, where ${\cal O}$ is a finite set of ordered pairs of
 non-negative integers (the {\it occurrence}
 in $FD$), $\nu$ is a function from $\cal O$ into $\Sigma$ ($\nu(i,j)$ is the
 value of $FD$ at occurrence $(i,j)$), and $p$ is a function from $\cal O$ into
 $\widetilde{h}=\{(r)a\rightarrow uf(r)|r\in J\}$ ($p(i,j)$ is the fuzzy
 production of $FD$ at occurrence $(i,j)$)
 satisfying the following conditions. There exists a
 sequence of words $(x_0,x_1,\cdots,x_m)$ in $\Sigma^*$ (called the trace of $FD$)
 such that $n\geq 1$ and

  (i) ${\cal O}=\{(i,j)|0\leq i\leq m\quad \mbox{and}\quad 1\leq j\leq |x_i|\}$,

  (ii) $\nu(i,j)$ is the $j$th symbol in $x_i$,

  (iii) For $0\leq i\leq m$, $x_{i+1}=u_1u_2\cdots u_{|x_i|}$, where $p(i,j)=
  r_{ij}\nu(i,j)\rightarrow u_jf(r_{ij})$ for $1\leq j\leq |x_i|$ and $r_{ij}$
  is determined by an ordinary production $\nu(i,j)\rightarrow u_j$. \
Then $FD$ is said to be a {\it fuzzy derivation} of $x_m$ from $x_0$, $m$ is
called the length of the fuzzy derivation $FD$.

\bigskip  If $x_0=\omega$, then $FD$ is said to be a fuzzy derivation of $x_m$ in $FG$.
  In general, there are more than one fuzzy derivations of $x_m$ in $FG$.

\bigskip  {\bf Definition 3.3}: The {\it grade of a string} $x$
 ($\in \Sigma^*$), denoted as $f_{FG}(x)$, is given as follows by using the
 concept of composition of fuzzy operations$^{[4]}$ and by the fuzzy
 derivations of $x$ in $FD$. Clearly, $f_{FG}(x)$
 is in $[0,1]$,
 $$f_{FG}(\omega)=\max_{r\in J}f(r),$$
 and for $x\neq\omega$,
 $$ f_{FG}(x)=\max\{\min f(r_{ij})|(i,j)\in {\cal O}\}$$
 where the maximum is taken over all fuzzy derivations of $x$ in $FG$.

\bigskip  {\bf Definition 3.4}: The {\it fuzzy language generated by fuzzy $0L$ system}
 $FG$ is
 $$\widetilde{L}(FG)=\{(f_{FG}(x),x)|x\in L(G)\}$$
 where $L(G)=\cup_{i\geq 0}h^i(\omega)$. Then $\widetilde{L}(G)$ is a fuzzy
 set over $\Sigma^*$.

\bigskip {\bf Definition 3.5}: A {\it fuzzy $E0L$ system} is a system
 $$FG=(\Sigma,\widetilde{h},\omega,J,f,\Delta)$$
 where $U(FG)=(\Sigma,\widetilde{h},\omega,J,f)$ is a fuzzy $0L$system and $\Delta\subseteq
 \Sigma$. The fuzzy language generated by the fuzzy $E0L$ system $FG$ is
 defined as
 $$\widetilde{L}(FG)=\{(f_{FG}(x),x)|x\in L(U(G))\cap \Delta^*\},$$
 where $L(U(G))=\cup_{i\geq 0}h^i(\omega)$ and $f_{FG}(x)=f_{U(FG)}(x)$.

\bigskip{\bf Definition 3.6}: A {\it fuzzy $T0L$ system} is a system
$$FG=(\Sigma,\widetilde{ H}, \omega,J,f),$$
where $\widetilde{H}$ is a nonempty finite set of fuzzy substitutions $\{\widetilde{h}_1,\cdots,\widetilde{h}_k\}$.
We take $\widetilde{H}=\cup_{i=1}^k\widetilde{h}_i$ and,
for every $\widetilde{h}_i$, $(\Sigma,\widetilde{h}_i,\omega,J_i,f_i)$
is a fuzzy $0L$ system, where $\{J_i\}_{i=1}^k$ are disjoint and $J=\cup_{i=1}^kJ_i$,
$f_i=f|_{J_i}$.

 By using the definition of derivation in a $T0L$ system ({\it cf.} p.~232
 of Ref. [6]), similar
 to Definition 3.2 and Definition 3.3, it is easy to define a fuzzy derivation
 in a fuzzy $T0L$ system and the grade $f_{FG}(x)$ of a string
 $ x(\in \Sigma^*)$.

Then we have:

\bigskip{\bf Definition 3.7}: The {\it fuzzy language generated by a fuzzy
 $T0L$ system} $FG$ is
$$\widetilde{L}(FG)=\{(f_{FG}(x),x)|x\in L(G)\},$$
where
\[
L(G)=\{x\in\Sigma^*|x=\omega\ \;\;\; {\rm or} \;\;\; x\in h_1\cdots h_m(\omega).
\]
In the above line $h_1,\cdots,h_m$
are substitutions corresponding to $\widetilde{h}_1,\cdots,\widetilde{h}_m\in
\widetilde{  H}$.

{\bf Definition 3.8}: A {\it fuzzy $ET0L$ system} is a system
$$FG=(\Sigma,\widetilde{H},\omega,J,f,\Delta)$$
where $U(FG)=(\Sigma,\widetilde{H},\omega,J,f)$ is a fuzzy $T0L$ system and $\Delta\subseteq \Sigma$.
The fuzzy language of $FG$ is defined by
$$\widetilde{L}(FG)=\{(f_{FG}(x),x)|x\in L(U(G))\cap \Delta^*\}$$
where
\[
L(U(G))=\{x\in\Sigma^*|x=\omega \;\;\; {\rm or} \;\;\; x\in h_1\cdots h_m(\omega)
\]
As before, $ h_1,\cdots,h_m$
are substitutions corresponding to $\widetilde{h}_1,\cdots,\widetilde{h}_m\in
\widetilde{  H}$.

  A fuzzy language is termed a fuzzy $0L$ (resp. $D0L$, $E0L$, $T0L$, $ET0L$)
  language if it coincides with 
 $\widetilde{L}(FG)$ for some fuzzy $0L$ (resp. $D0L$,
  $E0L$, $T0L$, $ET0L$) system $FG$.

  We use the notation ${\cal L}(F0L)$ (resp. ${\cal L}(FD0L)$, ${\cal L}(FE0L)$,
   ${\cal L}(FT0L)$, ${\cal L}(FET0L)$) to denote the set of all fuzzy $0L$
  (resp. $D0L$, $E0L$, $T0L$, $ET0L$) languages.

  From the above definitions the inclusion relations follow straightforwardly
  \begin{Proposition}
  $${\cal L}(FD0L)\subset{\cal L}(F0L)\subset{\cal L}(FE0L)
   \subset{\cal L}(FET0L),$$
  $${\cal L}(FD0L)\subset{\cal L}(F0L)\subset
   {\cal L}(FT0L)\subset{\cal L}(ET0L).$$
 \end{Proposition}

   Denoting by ${\cal L}(FCFL)$ (resp. ${\cal L}(FCSL)$) the
 set of all fuzzy context-free (resp. context sensitive) languages generated
 by some fuzzy grammar defined by Lee and Zadeh$^{[1]}$, we have
 \begin{Theorem} ${\cal L}(FCFL)\subset{\cal L}(FE0L)\subset{\cal L}(FET0L)
   \subset{\cal L}(FCSL)$.
 \end{Theorem}

  {\it Proof}. Let $\widetilde{L}=\widetilde{L}(FG)$ be a fuzzy context-free
  language which is generated by a fuzzy context-free grammar $FG=(V_N,V_T,
  P,S,J,f)$. Define a fuzzy $0L$ system
  $$FG'=(V_N\cup V_T,P\cup\{(r)a\rightarrow af_1(r) |a\in V_N\cup V_T\},S,J_1,f_1)$$
  where $J_1\backslash J$is the set of labels in $\{(r)a\rightarrow af_1(r)|a\in V_N\cup V_T\}$ and disjoint
  with $J$, $f_1:\quad (J_1\backslash J)\cup J\rightarrow [0,1]$ satisfies $f_1|_{J}=f$ and
  $f_1|_{J_1\backslash J}=\max_{r\in J}f(r)$.
  It is clear that $\widetilde{L}(FG)=\widetilde{L}(FG'')$, where
  $$FG''=(V_N\cup V_T,P\cup\{(r)a\rightarrow af_1(r) |a\in V_N\cup V_T\},S,J_1,f_1, V_T)$$
  is a fuzzy $E0L$ system. Hence each fuzzy context-free language is an fuzzy
  $E0L$ language.

  Let $FG_1=(\Sigma,\widetilde{H},\omega,J,f,\Delta)$ be an arbitrary fuzzy $ET0L$ system.
 Then $G_1=(\Sigma,H,\omega,$ $\Delta)$ is an ordinary $ET0L$ system,where $H=
 \{h|h\ \mbox{is the substitution corresponding to }\ \widetilde{h}\in\widetilde{H}\}$.
 From the proof of Theorem 10.5 in Herman and Rozenberg$^{[2]}$, there is a context
 sensitive grammar $G_2=(V_N,V_T,Q,\omega)$ such that $L(G_1)=L(G_2)$, where $L(G_1)$
 and $L(G_2)$ are the ordinary languages generated by $G_1$ and $G_2$ respectively.
 First we give an label to each production in $Q$. If a production $A\rightarrow \beta$
 in $Q$ has label $r$, then we define
 $$f'(r)=\min\{f(a)|a \ \mbox{is an symbol of }\ A\ \mbox{and rewritten by this production}\}.$$
 Then we can define a fuzzy context sensive grammar $FG_2=(V_N,V_T,Q,\omega,\{r\},f_1)$.
 It is clear that $f_{FG_1}(x)=f'_{FG_2}(x)$ for every $x\in L(G_1)=L(G_2)$.
 Hence $\widetilde{L}(FG_1)=\widetilde{L}(FG_2)$. Hence each fuzzy $ET0L$ language
 is a fuzzy context sensive language. $\Box$

\section{Formal language with the threshold $\lambda$}
{~~~~}Next we consider the relationship between fuzzy languages and ordinary languages.

\bigskip{\bf Definition 4.1}: Let $FG$ be a fuzzy $L$ system, $\lambda$ a real
number $0\leq\lambda< 1$, then a {\it language generated by $FG$ with the threshold}
$\lambda$ is defined by
$$L(FG,\lambda)=\{x\in\Sigma^*|f_{FG}(x)> \lambda\}.$$

  Then we have
 \begin{Theorem} For each fuzzy $T0L$ (resp. $ET0L$) system $FG$ with $f_i=c_i$
 ($i=1,\cdots,m$) being constants, then for any $\lambda$ ($0\le \lambda<1$),
 $L(FG,\lambda)$ is a $T0L$ (resp. $ET0L$) language.
 \end{Theorem}

  {\it Proof}. Let $H(\lambda)$ be the set of all substitutions which corresponds
  to $f_i>\lambda$ in $FG$. Since $FG$ satisfies $f_i=c_i$ ($i=1,\cdots,m$), from
  definitions, $L(FG,\lambda)$ is the language which is obtained from the only
  substitutions in $H(\lambda)$. It is obvious that our conclusion holds. $\Box$

 \section{Fuzzy entropy of a string}
{~~~~}In this section, we only consider fuzzy $0L$ systems and ordinary $E0L$ system.
 We continue
to view a substitution $h$ as a set of productions for each symbol in the alphabet.
For a string $w$ and a symbol $b$, $N_b(w)$ denotes the number of occurrences of
symbol $b$ in $w$. Let $FG=(\Sigma,\widetilde{h}, \omega,J,f)$ be a fuzzy $0L$ system, then
for each $b$ in $\Sigma$, $R_b$ means the set of labels of the rules having $b$
on the left-hand side.
For $w=a_1a_2\cdots a_m$ in $\Sigma^*$, $A_w$ means a subset of $J^*$
defined as follows:
$$A_w=\{\alpha=\pi_1\cdots \pi_m|(\pi_i)a_i\rightarrow u_if(\pi_i)\ \mbox{in}\ \widetilde{h},1\le i\le m\}.$$
For a set $F$, \#$F$ denotes the cardinality of $F$. Define $d(a)=\#R_a$ which
is called the {it degree} of $a$. For any $a\in \Sigma$, $r\in R_a$, we define
$$\mu(r)=\frac{f(r)}{\sum_{r_i\in R_a}f(r_i)}.$$ Then it is obvious that $\sum_{r\in R_a}\mu(r)=1$.
If for any $r_1,r_2\in R_a$, $f(r_1)=f(r_2)$, then we have $\mu(r)=\frac{1}{d(a)}$,
for any $r\in R_a$. If $\mu(r)=\frac{1}{d(a)}$ for all $a\in \Sigma,r\in R_a$
holds, then $FG$ is called {\it uniformly fuzzy $0L$ system}.
 For a $L$ system $G$, $L(G)$ denote the
ordinary language generated by $G$.

  \bigskip{\bf Definition 5.1}: Let $FG=(\Sigma,\widetilde{h},\omega,J,f)$ be a fuzzy $0L$
 system, $E_f$ is defined as follows: for any $w$ in $\Sigma^*$,
 $$E_f(\Lambda)=0,\quad E_f(w)=-\sum_{\alpha\in A_w}\mu(\alpha)\log_2\mu(\alpha),$$
 where $f(\pi\alpha)=\mu(\pi)\mu(\alpha)$, for $\alpha\in J^*$ and
 $\pi\in J$. $E_f(w)$ is termed the {\it fuzzy entropy} of $w$.

  For some $c\geq 0$, we define $L(FG,E_f\le c)=\{w\in\Sigma^*|w\in\cup_{i\geq 0}h^i(\omega)\
  \mbox{and}\ E_f(w)\leq c\}$. Then $L(FG,E_f\le c)$ is called {\it bounded fuzzy
  entropy language}. Denote ${\cal L}(BFEL)$ the class of languages with bounded
  fuzzy entropy. We use ${\cal L}(RGL)$ and ${\cal L}(E0L)$ denote the class of
  regular languages and of languages generated by an $E0L$ system respectively.

\begin{Theorem} For any $E0L$ system $G$, there exists a uniformly fuzzy $0L$ system $FG$
 such that $L(G)=L(FG,E_f\le 0)$. Hence
${\cal L}(E0L)\subseteq {\cal L}(BFEL)$.
\end{Theorem}
  {\it Proof}. By Lemma 1 of Ref.[12], let $G'=(\Sigma',h',\omega',\Delta')$
 be the $E0L$ system such that (1)for all $a\in \Delta'$, $d(a)=1$ and for all
 $a\in \Sigma'\backslash\Delta'$, $d(a)\geq 2$, and (2) $L(G)=L(G')$ hold.
 Consider a uniformly fuzzy $0L$ system
 $FG=(\Sigma',\widetilde{h'},\omega',J,f)$. From the definition of language generated by
 a $E0L$ system, it is sufficient to show that $E_f(w)=0$ if and only if $w\in \Delta'^*$.
 Let $\Sigma'=\{a_1,\cdots,a_k\}$ and $\Delta'=\{a_1,\cdots,a_r\}$, then
 $V_s=(\log_2 d(a_1),\cdots,\log_2 d(a_k))^T=(0,\cdots,0,\nu_{r+1},\cdots,\nu_k)^T$,
 where $\nu_i\ge 1$ ($r+1\le i\le k$). Let $U_s=(N_{a_1}(w),  \cdots, N_{a_k}(w))$.
 Since $FG$ is uniformly fuzzy $0L$ system, it is easy to see that
 \begin{eqnarray*}
& & E_f(w)=0 \Longleftrightarrow U_sV_s=0\\
& &\Longleftrightarrow N_{a_{r+1}}(w)\nu_{r+1}+\cdots+N_{a_k}(w)\nu_k=0\\
& & \Longleftrightarrow N_{a_i}(w)=0,\quad (\mbox{for}\quad i=r+1,\cdots,k)
 \Longleftrightarrow w\in \Delta'^*.
 \end{eqnarray*}
 Then we obtain
 $$L(G)=L(G')=L(FG,E_f\le 0).$$
 $\Box$

     Now we are in a position to propose:
  \begin{Theorem} ${\cal L}(BFEL)= {\cal L}(E0L)$.
\end{Theorem}
  {\it Proof}. From Theorem 3, it is sufficient to prove ${\cal L}(BFEL)
  \subseteq {\cal L}(E0L)$.
  Let $FG=(\Sigma,\widetilde{h},\omega,J,f)$ be a fuzzy $0L$ system. Denote $S=(\Sigma,
  h,\omega, \mu)$, then $S$ is a stochastic $0L$ system defined in Ref.[10].
  For a given constant $c$ ($\geq 0$),
  it is easy to see $L(FG,E_f\le c)=L(S,c)$, where $L(S,c)$ is also defined
  in Ref.[10]. Then from the Proof of Theorem of Ref.[10], we have $L(FG,E_f\le c)$
  is in ${\cal L}(E0L)$.
   $\Box$

\bigskip  {\bf Remark}. Using the same method, one can introduce fuzziness into
  many other families of $L$ system.

\section*{ACKNOWLEDGMENTS}

{~~~~}The author would like to express his thanks to Prof. Bai-lin Hao for
useful discussions and reading the manuscript, and to Prof. Wei-mou Zheng,
Dr. Hua-lin Shi,
Dr. Shuang-hu Wang and Dr. Guo-yi Cheng for many suggestions.


\begin{thebibliography}{99}
\renewcommand{\baselinestretch}{0.2}
{\small
\bibitem[1]{} N. Chomsky, Syntactic Structures (Mouton The Hague, 1957).

\bibitem[2]{} G.T. Herman and G. Rozenberg, Developmental Systems and Languages
(North-Holland Publ. Amsterdam, 1975).

\bibitem[3]{} E.T. Lee and L.A. Zadeh, Note on fuzzy languages, 
 Information Science {\bf  1} (1969) 421-434.

\bibitem[4]{} A. Lindenmayer,  Mathematical models for cellular interactions
in development, I \& II,  J. Theor. Biol. {\bf 18} (1968) 280-315.

\bibitem[5]{} M. Mizumoto, J. Toyoda and K. Tanaka, $N$-fold fuzzy grammars,
Information Science {\bf 5} (1973) 25-43.

\bibitem[6]{} G. Rozenberg and A. Salomaa, The Mathematical Theory of $L$ Systems
(Academic Press, New York, 1980).

\bibitem[7]{} A. Salomaa, Formal Languages (Academic Press, New York/London,
1973).

\bibitem[8]{} A. Salomaa, Jewels of Formal Language Theory (Computer Science
Press, Rockville, 1981).

\bibitem[9]{} Hui-min Xie, Grammatical Complexity and One-dimensional
Dynamical Systems (World Scientific, 1996).

\bibitem[10]{} T. Yokomo, Stochastic characterizations of $E0L$ languages,
  Information and Control {\bf 45} (1980) 26-33.

\bibitem[11]{} L.A. Zadeh, Fuzzy sets,  Information and Control {\bf 8}
(1965) 338-353.

\bibitem[12]{} H. -J. Zimmermann, Fuzzy Set Theory and its Applications (Boston:
Kluwer-Nijhoff Publ. 1985).

}
\end{thebibliography}
\end{document}